\documentstyle[12pt,leqno]{article}
\def\date{\today}
\newtheorem{theorem}{Theorem}[section]
\newtheorem{definition}[theorem]{Definition}

\newtheorem{lemma}[theorem]{Lemma}

\newtheorem{proposition}[theorem]{Proposition}

\catcode`\@=11
\def\@begintheorem#1#2{\it \trivlist \item[\hskip
 \labelsep{\sc #1\ #2.}]}
\def\@opargbegintheorem#1#2#3{\it \trivlist\item[\hskip%
 \labelsep{\sc #1\ #2.\ (#3)}]}
\def\@endtheorem{\endtrivlist}

\def\@listI{\leftmargin\leftmargini \parsep 1pt plus 2.5pt
 minus 1pt\topsep 5pt plus 2pt minus 3pt%
 \itemsep 0pt plus 2.5pt minus 1pt}
\let\@listi\@listI
\@listi

\def\@sect#1#2#3#4#5#6[#7]#8{\ifnum #2>\c@secnumdepth%
 \def \@svsec {}\else \refstepcounter {#1}\edef \@svsec%
 {\csname the#1\endcsname. \hskip .1em }\fi \@tempskipa%
 #5\relax \ifdim \@tempskipa >\z@ \begingroup #6\relax%
 \@hangfrom {\hskip #3\relax \@svsec }{\interlinepenalty%
 \@M #8\par }\endgroup \csname #1mark\endcsname {#7}%
 \addcontentsline {toc}{#1}{\ifnum #2>\c@secnumdepth%
 \else \protect \numberline {\csname the#1\endcsname. }%
 \fi #7}\else \def \@svsechd {#6\hskip #3\@svsec #8%
 \csname #1mark\endcsname {#7}\addcontentsline {toc}{#1}%
 {\ifnum #2>\c@secnumdepth \else \protect \numberline%
 {\csname the#1\endcsname. }\fi #7}}\fi \@xsect {#5}}

\def\section{\@startsection {section}{1}{\z@ }%
 {-3.5ex plus -1ex minus -.2ex}{2.3ex plus .2ex}{\sc }}

\def\thebibliography#1{%
 \section *{References\@mkboth {REFERENCES}{REFERENCES}}%
 \list {[\arabic {enumi}]}{\settowidth \labelwidth {[#1]}%
 \leftmargin \labelwidth \advance \leftmargin \labelsep %
 \usecounter {enumi}} \def \newblock %
 {\hskip .11em plus .33em minus -.07em} \sloppy \clubpenalty 4000%
 \widowpenalty 4000 \sfcode`\.=1000\relax}

\def\@maketitle{%
 \newpage \null \vskip 2em
 \begin{center}{\Large\sc \@title \par }
 \vskip 1.5em
 {\large \lineskip .5em
 \begin {tabular}[t]{c}\@author
 \end{tabular}\par }
 \end{center}\par \vskip 1.5em}

\def\abstract{%
\if@twocolumn \section *{Abstract}
 \else \small\quotation\noindent{\sc Abstract.}\fi}

\catcode`\@=13

\def\pf{\noindent{\sc Proof. }}
\def\qed{\hspace*{\fill}
\mbox{\hphantom{mm}{\sc q.e.d.}}\\
\par\vskip-5mm}

\def\MU{\hskip.2mm\cdot\hskip.8mm}
\def\boldk{{\bf k}} \def\bk{{\bf k}}
\def\tp#1#2{{#1}^{\otimes{#2}}}
\def\ident{1\!\!1} \def\id{\ident}
\def\D{\Delta}
\def\ot{\otimes}
\def\oV{{\overline V}}
\def\Ker{{\mbox{\rm Ker}}}
\def\br#1{{\mbox{\rm Br}_{#1}}}
\def\b{*}
\def\Sum{\mbox{$\sum$\hskip1mm}}
\def\bigodo{\mbox{$\bigodot$}}
\def\od{\odot}
\def\bod{\bigodo}
\def\bv{{\cal B}_\b(V)}
\def\db{d_{\cal B}}
\def\calb{{\cal B}}
\def\bodvbv{\bod(V,\bv)}
\def\susp{\uparrow\kern -.05em}
\def\Der{\mbox{\rm Der}}
\def\der#1{\mbox{\rm Der}^{#1}_V(\bod(V,\bv))}
\def\desusp{\downarrow\kern -.05em}
\def\calr{{\cal R}}
\def\dr{d_{\cal R}}
\def\cals{{\cal S}}
\def\ds{d_{\cal S}}
\def\Im{\mbox{\rm Im}} \def\bolds{{\bf s}}
\def\L{{\cal L}} \def\boldt{{\bf t}}
\def\cb#1{\calb_{#1}(V)} 
\def\cc#1#2{CC_{#1}(K_{#2})}
\def\Br#1#2{\mbox{\rm Br}_{#2}(#1)}
\def\hh#1#2#3{(\id^{\odot(#1)}\odot {#2}_N\odot \id^{\odot(#3)})}
\def\e{{\bf e}}

\begin{document}
\baselineskip18pt

\title{COHOMOLOGY OF DRINFEL'D ALGEBRAS:
\\ A GENERAL NONSENSE APPROACH}
\author{\rm Martin MARKL%
\thanks{Partially supported by the grant AV \v CR \#119 105 and
by the Gelbart Research Institute, Bar Ilan University, \hskip7mm
Israel}
\hskip1mm and Steve SHNIDER}

\maketitle

\section{Preliminaries}

Recall that a {\em Drinfel'd algebra\/} (or a
{\em quasi-bialgebra\/} in the
original terminology of~\cite{drinfeld:Alg.iAnaliz89}) is an object
$A=(V,\MU,\Delta,\Phi)$, where
$(V,\MU,\Delta)$ is an associative, not necessarily coassociative,
unital and counital $\boldk$-bialgebra, $\Phi$ is an
invertible element of
$\tp V3$,
and the usual coassociativity property is replaced by
the condition which we shall refer to as quasi-coassociativity:
\begin{equation}
\label{14}
(\ident\otimes\Delta)\Delta \cdot \Phi
= \Phi\cdot(\Delta\otimes\ident)\Delta,
\end{equation}
where we use the dot $\MU$ to indicate both the (associative)
multiplication on $V$ and the induced multiplication
on $V^{\otimes
3}$.
Moreover, the validity of the following
``pentagon identity'' is required:
\[
(\ident^2\otimes\Delta)(\Phi)\cdot(\Delta\otimes\ident^2)
(\Phi) =
(1\otimes\Phi)\cdot(\ident
\otimes\Delta\otimes\ident)(\Phi)\cdot(\Phi\otimes1),
\]
$1\in V$ being the unit element and $\ident$, the identity map
on $V$. If $\epsilon :V \to \bk$ ($\bk$ being the ground field) is the
counit of the coalgebra $(V,\D)$ then, by definition, $(\epsilon \ot
\id)\Delta = (\id \ot \epsilon)\D = \id$.
We have a natural
splitting $V = \oV \oplus \bk$, $\oV := \Ker(\epsilon)$,
given by the embedding $\bk \to V$, $\bk\ni c \mapsto c\cdot 1 \in V$.

For a $(V,\MU)$-bimodule $N$, recall the following generalization of
the $M$-{\em construction\/}
of~\cite[par.~3]{markl-stasheff:Journ.ofAlgToApp} introduced
in~\cite{markl-shnider:preprint}.
Let $F^\b=\bigoplus_{n\geq 0}F^n$ be
the free unitary nonassociative $\boldk$-algebra generated by
$N$, graded by the length of words. The space $F^n$
is the direct sum of
copies of $N^{\ot n}$ over the set $\br n$ of full bracketings of $n$
symbols, $F^n = \bigoplus_{u\in {\rm Br}_n}N^{\ot n}_u$.
For example, $F^0=\bk$, $F^1 = N$, $F^2 =
N^{\ot 2}$, $F^3 = N^{\ot 3}_{(\bullet\bullet)\bullet} \oplus
N^{\ot 3}_{\bullet(\bullet\bullet)}$,~etc.
The algebra $F^\b$ admits a
natural left action,
$(a,f)\mapsto a\bullet f$, of the algebra $(V,\MU)$ given by the
rules:
\begin{itemize}
\item[(i)]
on $F^0=\bk$, the action is given by the augmentation $\epsilon$,
\item[(ii)]
on $F^1=N$, the action is given by the left action of $V$ on $N$
and
\item[(iii)]
$a\bullet (f\star g)= \sum(\Delta'(a)\bullet
f)\star(\Delta''(a)\bullet g)$,
\end{itemize}
where $\star$ stands for the multiplication in $F^\b$ and we use the
Sweedler
notation $\Delta(a)= \sum \Delta'(a)\otimes \Delta''(a)$. The right
action $(f,b)\mapsto f\bullet b$ is defined by similar rules. These
actions define on $F^\b$ the structure of a $(V,\MU)$-bimodule.

Let $\sim$ be the relation on $F^\b$ $\star$-multiplicatively
generated
by the expressions of the form
\[
\sum\left((\Phi_1\bullet x)\star\left((\Phi_2\bullet
y)\rule{0mm}{4mm}\star(\Phi_3\bullet
z)\right)\right)
\sim\sum\left(\left((x\bullet\Phi_1)\star\rule{0mm}{4mm}
(y\bullet\Phi_2)\right)\star
(z\bullet\Phi_3)\right),
\]
where $\Phi = \sum\Phi_1\otimes \Phi_2\otimes\Phi_3$ and $x,y,z\in
F^\b$.
Put $\bigodo(N):= F/\sim$.
Just as in~\cite[Proposition~3.2]{markl-stasheff:Journ.ofAlgToApp}
one proves that the $\bullet$-action induces on $\bigodo(N)$
the structure of
a $(V,\MU)$-bimodule (denoted again by~$\bullet$) and that
$\star$ induces
on $\bigodo(N)$ a nonassociative multiplication denoted by
$\odot$. The operations are related by
\[
a\bullet (f\od g)= \Sum(\Delta'(a)\bullet f)\od(\Delta''(a)\bullet g)
\mbox{ and }
(f\od g)\bullet b = \Sum(f\bullet \Delta'(b))\od (g \bullet
\Delta''(b)),
\]
for $a,b\in V$ and $f,g\in \bigodo(N)$.
The multiplication $\odot$ is quasi-associative in the sense
\begin{equation}
\label{3}
\sum(\Phi_1\bullet x)\odot\left((\Phi_2\bullet
y)\rule{0mm}{4mm}\odot(\Phi_3\bullet
z)\right)
=\sum\left((x\bullet\Phi_1)\odot\rule{0mm}{4mm}
(y\bullet\Phi_2)\right)\odot
(z\bullet\Phi_3)
\end{equation}
The construction described above is functorial in the sense that any
$(V\MU)$-bimodule map $f:N' \to N''$ induces a natural
$(V\MU)$-linear algebra homomorphism $\bod(f):\bod(N')\to \bod(N'')$.
As it was shown in~\cite{markl-stasheff:Journ.ofAlgToApp},
for any $(V,\MU)$-bimodule $N$
there exists a natural homomorphism of
$\bk$-modules $J= J(N): \bod(N)\to \bigotimes (N)$. If $f:N' \to N''$
is as above then $J(N'')\circ \bod(f)= \bigotimes(f)\circ J(N')$.

Since the defining relations~(\ref{3}) are homogeneous with
respect to length,
the grading of $F^\b$ induces on $\bod(N)$ the grading $\bod^\b(N)=
\bigoplus_{i\geq 0}\bod^i(N)$. If $N$ itself is a graded vector
space, we have also the obvious second grading, $\bod(N) =
\bigoplus_{j}\bod(N)^j$, which coincides with the first grading if
$N$ is concentrated in degree~1.

Let $\mbox{\rm Der}^n_V(\bod(N))$ denote the set of
$(V,\MU)$-linear derivations of degree $n$ (relative to the
second grading) of the (nonassociative) graded algebra
$\bod(N)^\b$. One sees immediately that there is an one-to-one
correspondence between the elements $\theta \in\mbox{\rm
Der}^n_V(\bod N)$
and $(V,\MU)$-linear homogeneous degree~$n$ maps $f:N^\b\to
\bod(N)^\b$.

If $N=X\oplus Y$, then $\bod(X\oplus Y)$ is naturally
bigraded, $\bod^{\b,\b}(X\oplus Y)= \bigoplus_{i,j\geq
0}\bod^{i,j}(X\oplus Y)$, the bigrading being defined by saying that
a monomial $w$ belongs to $\bod^{i,j}(X\oplus Y)$ if there are exactly
$i$
(resp.~$j$) occurrences of the elements of $X$ (resp.~$Y$) in $w$.
If $X,Y$ are graded vector spaces then there is a second bigrading
$\bod(X\oplus Y)^{\b,\b}= \bigoplus_{i,j}
\bod(X\oplus Y)^{i,j}$ just as above.

Let $(\bv,\db)$ be the (two-sided) normalized
bar resolution of the algebra
$(V,\MU)$ (see~\cite[Chapter~X]{maclane:homology}),
but considered with the opposite grading.
This means that $\bv$ is
the graded $(V\MU)$-bimodule, $\bv = \bigoplus_{n\leq
1}\calb_n(V)$, where $\calb_1(V):=V$ with the $(V,\MU)$-bimodule
structure induced by the multiplication $\MU$,
$\calb_0(V):=V\otimes V$ (the free $(V,\MU)$-bimodule on $\bk$), and
for $n\leq- 1$, $\calb_n(V)$ is the free $(V,\MU)$-bimodule on
$\oV^{\otimes( -n)}$,
i.e.~the vector space $V\ot \oV^{\otimes(-n)}\ot V$
with the action of ($V,\MU)$ given by
\[
u\cdot (a_0\otimes\cdots\otimes a_{-n+1}):= (u\cdot
a_0\otimes\cdots\otimes a_{-n+1})\
\mbox{ and }\
(a_0\otimes\cdots\otimes a_{-n+1})\cdot w := (a_0\otimes\cdots\otimes
a_{-n+1}\cdot w)
\]
for $u,v,a_0,a_{-n+1}\in V$ and $a_1,\ldots,a_{-n}\in \oV$.
If we use the more compact
notation (though a nonstandard one), writing $(a_0|\cdots|a_{-n+1})$
instead of $a_0\otimes\cdots\otimes a_{-n+1}$, the differential $\db
: \calb_n(V)
\to \calb_{n+1}(V)$ is, for $n\leq 0$, defined as
\[
\db(a_0|\cdots|a_{-n+1}):=
\sum_{0\leq i\leq -n}(-1)^{i}(a_0|\cdots|a_i\cdot a_{i+1}|
\cdots|a_{-n+1}).
\]
Here, as is usual in this context, we make no distinction between
the elements of $V/\bk \cdot 1$ and their representatives in $\oV$. We
use the same convention throughout all the paper.
Notice that the differential $\db$ is a $(V,\MU)$-bimodule map.

Put $\bodvbv := \bod(\susp V \oplus \susp \bv)$, where $\susp$
denotes, as usual, the suspension of a graded vector space and $V$ is
interpreted as a graded vector space concentrated in degree zero.
Let $\der i$
denote, for each $i$, the space of degree $i$ derivations of the
algebra $\bodvbv$ which are also $(V,\MU)$-linear maps.
Let us define the derivation $D_{-1}\in \der 1$ by
$D_{-1}|_{\susp \bv}:= \susp \db \desusp$ and
$D_{-1}|_{\susp V}:= 0$.
Clearly $D_{-1}(\bod(V,{\cal B}_\b(V))^{i,j})
\subset \bod(V,{\cal B}_\b(V))^{i,j+1}$ and
$D_{-1}(\bod^{i,j}(V,\bv))\subset \bod^{i,j}(V,\bv)$
for any $i,j \geq 0$. We also see immediately that
$D_{-1}^2=0$.

Let us consider, for any $n\geq 1$,
the complex $(\bod^{n-1,1}(V,\calb_\b(V)),D_{-1})$,
i.e.~the complex
\[
0\longleftarrow \bod^{n-1,1}(V,V)
\stackrel{D_{-1}}{\longleftarrow} \bod^{n-1,1}(V,\calb_0(V))
\stackrel{D_{-1}}{\longleftarrow}\bod^{n-1,1}(V,\calb_{-1}(V))
\longleftarrow\cdots
\]
\begin{lemma}
The complex $(\bod^{n-1,1}(V,\calb_\b(V)),D_{-1})$ is acyclic, for
any $n\geq 1$.
\end{lemma}

\pf We have the decomposition \[
\bod^{n-1,1}(V,\calb_\b(V)) =
\bigoplus_{1\leq i\leq n}\bod^{n-1,1}_i(V,\calb_\b(V)),
\]
where $\bod^{n-1,1}_i(V,\calb_\b(V))$ denotes the subspace of
$\bod^{n-1,1}(V,\calb_\b(V))$ spanned by monomials having an
element of $\calb_\b(V)$ at the $i$-th place. The differential
$D_{-1}$
obviously respects this decomposition and the canonical isomorphism
$J$ of~\cite{markl-stasheff:Journ.ofAlgToApp} mentioned above
identifies $\bod^{n-1,1}_i(V,\calb_\b(V))$ to
$\tp V{(i-1)}\ot \bv \ot \tp V{(n-1)}$. Under this identification the
differential $D_{-1}$ coincides with
$\tp \id{(i-1)}\ot \db\ot \tp\id{(n-i)}$ and the rest follows from the
K\"unneth formula and the acyclicity of $(\bv,\db)$.
\qed

\section{Properties of $\der \b$}

Let $C = (C,\MU,1_C)$ be a unital associative algebra and let $C
\stackrel{\epsilon}{\longleftarrow}(\calr,\dr)$, $(\calr,\dr)=
R_0 \stackrel{\dr}{\longleftarrow}R_1
\stackrel{\dr}{\longleftarrow}\cdots$, be a complex of free
$C$-bimodules (we consider $C$ as a $C$-bimodule with the bimodule
structure induced by the multiplication). Similarly, let $D
\stackrel{\eta}{\longleftarrow}(\cals,\ds)$ with $(\cals,\ds)=
S_0 \stackrel{\ds}{\longleftarrow}S_1
\stackrel{\ds}{\longleftarrow}\cdots$, be an acyclic complex of
$C$-bimodules. To simplify the notation, we write sometimes
$R_{-1}$ (resp.~$S_{-1}$, resp.~$\dr$, resp.~$\ds$) instead of $C$
(resp.~$D$, resp.~$\epsilon$, resp.~$\eta$).
Let
\[
Z := \{f=(f_i)_{i\geq -1};\ f_i:R_i\to S_i\ \mbox{a $C$-bimodule
map and}\ f_i\circ \dr = \ds \circ f_{i+1}\
\mbox{for any $i\geq -1$}\}.
\]
Let us define, for a sequence $\chi = (\chi_i)_{i\geq -1}$ of
$C$-bimodule maps $\chi_i: R_i \to S_{i+1}$,
$\nabla(\chi)= (\nabla(\chi)_i)_{i\geq -1}\in Z$ by
$\nabla(\chi)_i:= \ds \circ \chi_i+\chi_{i-1}\circ\dr$. Let $B:=
\Im(\nabla)\subset Z$. For a $C$-bimodule $M$ let $M_I$ denote the set
of invariant elements of $M$, $M_I:= \{x\in M;\ cx = xc\ \mbox{for
any $c\in C$}\}$.

\begin{lemma}
\label{1}
Under the notation above, the correspondence
$Z\ni f=(f_i)_{i\geq -1}\mapsto f_{-1}(1_C)\in D_I$ induces an
isomorphism $\Omega : Z/B\cong D_I/\eta({S_0}_I)$. Moreover, if
$f_{-1}(1_C)=\eta(h)$ for some $h\in {S_0}_I$ then $f=\nabla(\chi)$
for some $\chi = (\chi_i)_{i\geq -1}$ with $\chi_{-1}(1)=h$.
\end{lemma}

\pf We show first that $\Omega$ is well-defined. If
$f=\nabla(\chi)$ then $f_{-1}=
\eta\circ \chi_{-1}$, therefore $f_{-1}(1_C)= \eta(h)$ with $h:=
\chi_{-1}(1_C)\in {S_0}_I$ and $\Omega(f)=0$.

Let us prove that $\Omega$ is an epimorphism. For $z\in D_I$
define a $C$-bimodule map $f_{-1}: C\to D$ by $f_{-1}(c):= cz$ (=
$zc$)
for $c\in C$. Because $(\calr,\dr)$ is free and $(\cals,\ds)$ is
acyclic,
$f_{-1}$ lifts to some $f=(f_i)_{i\geq -1}\in Z$ by
standard homological
arguments~\cite[Theorem~III.6.1]{maclane:homology}.

It remains to prove that $\Omega$ is a monomorphism. For
$f=(f_i)_{i\geq -1} \in Z$, $\Omega(f)=0$ means that $f_{-1}(1_C)=
\eta(h)$ for some $h\in {S_0}_I$. The $C$-bimodule map $\chi_{-1}:
C\to S_0$ defined by $\chi_{-1}(c):= ch$ (= $hc$) for $c\in C$ clearly
satisfies $f_{-1}= \eta \circ \chi_{-1}$. A standard homological
argument (see again \cite[Theorem~III.6.1]{maclane:homology}
then enables one to extend $\chi_{-1}$ to a `contracting
homotopy' $\chi = (\chi_i)_{i\geq -1}$ with $f= \nabla(\chi)$.
\qed

\begin{definition}
\label{0}
For $n\geq 2$ and $k\geq 0$ let $J_k(n)$ be the subspace of
$\der {n-1-k}$ consisting of derivations $\theta$ satisfying
\begin{itemize}
\item[(i)]
$\theta(\bod(V,\bv)^{i,j})\subset \bod(V,\bv)^{i+n-1,j-k}$,
\item[(ii)]
$\theta(\bod^{i,j}(V,\bv))\subset \bod^{i+n-1,j}(V,\bv)$,
\item[(iii)]
$[D_{-1},\theta]=0$ if $k=0$ and $\theta|_{\susp \calb_1(V)}=0$ if
$k\geq 1$.
\end{itemize}
\end{definition}

Let us observe that, for $\theta \in J_{\geq 1}(n)$, $\theta|_{\susp
V}=0$.
This follows from item~(i) of the definition above.
Observe also that $J_\b(n)$ is
invariant under the differential $\nabla$ defined by
$\nabla(\theta):= [D_{-1},\theta]$,
$\nabla(J_k(n))\subset J_{k-1}(n)$ for $k\geq 1$ and
$\nabla(J_0(n))=0$.

\begin{proposition}
\label{12}
$H_{\geq 1}(J_\b(n),\nabla)=0$ while
\begin{equation}
\label{2}
H_0(J_\b(n),\nabla) =
\bod^{n-1,1}(V,\calb_1(V))_I
\oplus \bod^n(V)_I.
\end{equation}
\end{proposition}

\pf
Let $k>0$ and let $\theta \in J_k(n)$. As $\theta|_{\susp V}=0$,
$\theta$ is
given by its restriction to $\bv$, namely by a sequence of
$(V,\MU)$-bimodule maps
$\theta_i:\calb_i(V)\to \bod^{n-1,1}(V,\calb_{i-k}(V))$, $i\leq 1$.
Suppose that $\theta$ is a $\nabla$-cocycle, i.e.~that
$\nabla(\theta)=0$. This means that the diagram
\[
\unitlength=1.45mm
\begin{picture}(100.00,25.00)(25,10)
\put(40.00,10.00){\makebox(0,0)[cc]{$\bod^{n-1,1}(V,\calb_{1-k}(V))$}}
\put(70.00,10.00){\makebox(0,0)[cc]{$\bod^{n-1,1}(V,\calb_{-k}(V))$}}
\put(100.00,10.00){\makebox(0,0)[cc]{$\bod^{n-1,1}(V,%
\calb_{-1-k}(V))$}}
\put(40.00,30.00){\makebox(0,0)[cc]{$\calb_1(V)=V$}}
\put(70.00,30.00){\makebox(0,0)[cc]{$\calb_0(V)$}}
\put(100.00,30.00){\makebox(0,0)[cc]{$\calb_{-1}(V)$}}
\put(40.00,26.00){\vector(0,-1){11.00}}
\put(70.00,26.00){\vector(0,-1){11.00}}
\put(55.00,35.00){\makebox(0,0)[cc]{$\db$}}
\put(85.00,35.00){\makebox(0,0)[cc]{$\db$}}
\put(55.00,15.00){\makebox(0,0)[cc]{$D_{-1}$}}
\put(85.00,15.00){\makebox(0,0)[cc]{$D_{-1}$}}
\put(118.00,15.00){\makebox(0,0)[cc]{$D_{-1}$}}
\put(45.00,21.00){\makebox(0,0)[cc]{$\theta_1$}}
\put(76.00,21.00){\makebox(0,0)[cc]{$\theta_0$}}
\put(105.00,21.00){\makebox(0,0)[cc]{$\theta_{-1}$}}
\put(20.00,30.00){\makebox(0,0)[cc]{0}}
\put(20.00,10.00){\makebox(0,0)[cc]{0}}
\put(100.00,26.00){\vector(0,-1){10.89}}
\put(115.00,35.00){\makebox(0,0)[cc]{$\db$}}
\put(128.00,30.00){\makebox(0,0)[cc]{$\cdots$}}
\put(128.00,10.00){\makebox(0,0)[cc]{$\cdots$}}
\put(30.00,30.00){\vector(-1,0){8.00}}
\put(64.00,30.00){\vector(-1,0){15.93}}
\put(94.00,30.00){\vector(-1,0){17.93}}
\put(124.00,30.00){\vector(-1,0){17.93}}
\put(26.00,10.00){\vector(-1,0){4.00}}
\put(86.00,10.00){\vector(-1,0){4.00}}
\put(57.00,10.00){\vector(-1,0){4.1}}
\put(124.00,10.00){\vector(-1,0){9}}
\end{picture}
\]
is commutative. Since $\theta_1=0$ by item~(iii) of
Definition~\ref{0}, Lemma~\ref{1} (with $C = (V,\MU,1_V)$ and
$D=\bod^{n-1,1}(V,\calb_{1-k}(V))$) gives a sequence
$\chi_i:\calb_i(V)\to \bod^{n-1,1}(V,\cb {i-k-1})$ of
$(V,\MU)$-bimodule maps, $i\leq 1$, such that $\theta_i=
D_{-1}\circ\chi_i + \chi_{i+1}\circ \db$. We can, moreover, suppose
that $\chi_1=0$, thus
the sequence $(\chi_i)_{i\leq 1}$ determines a derivation
$\chi \in J_{k+1}(n)$ with $\nabla(\chi)=\theta$. This proves
$H_k(J_\b(n),\nabla)=0$ for $k\geq 1$.

A derivation $\theta \in J_0(n)$ is given by two independent data: by
the restriction $\theta_V := \theta|_{\susp V}: V\to\bod^n(V)$ and by
the restriction
$\theta_{\bv}:= \theta|_{\susp \bv}: \bv \to \bod^{n-1,1}(V,\bv)$.
As $D_{-1}|_{\susp V}=0$, the
condition $\nabla(\theta)= [D_{-1},\theta]=0$
imposes no restrictions
on $\theta_V$ and, because $\nabla(\chi)|_{\susp V}=0$ for any
$\chi \in J_1(n)$, the contribution of $\theta_V$ to
$H_0(J_\b(n),\nabla)$ is parametrized by $\theta_V(1_V)$, i.e.~by an
element of $\bod^n(V)_I$. This explains the second summand
in~(\ref{2}).

The restriction $\theta_{\bv}$ is in fact a sequence $\theta_i:\cb i
\to \bod^{n-1,1}(V,\cb i)$, $i\leq 1$, of $(V,\MU)$-bimodule maps and
the condition $\nabla(\theta)=0$ means that the diagram
\[
\unitlength=1.45mm
\begin{picture}(100.00,25.00)(25,10)
\put(40.00,10.00){\makebox(0,0)[cc]{$\bod^{n-1,1}(V,\calb_{1}(V))$}}
\put(70.00,10.00){\makebox(0,0)[cc]{$\bod^{n-1,1}(V,\calb_{0}(V))$}}
\put(100.00,10.00){\makebox(0,0)[cc]{$\bod^{n-1,1}(V,\calb_{-1}(V))$}}
\put(40.00,30.00){\makebox(0,0)[cc]{$\calb_1(V)=V$}}
\put(70.00,30.00){\makebox(0,0)[cc]{$\calb_0(V)$}}
\put(100.00,30.00){\makebox(0,0)[cc]{$\calb_{-1}(V)$}}
\put(40.00,26.00){\vector(0,-1){11.00}}
\put(70.00,26.00){\vector(0,-1){11.00}}
\put(55.00,35.00){\makebox(0,0)[cc]{$\db$}}
\put(85.00,35.00){\makebox(0,0)[cc]{$\db$}}
\put(55.00,15.00){\makebox(0,0)[cc]{$D_{-1}$}}
\put(85.00,15.00){\makebox(0,0)[cc]{$D_{-1}$}}
\put(118.00,15.00){\makebox(0,0)[cc]{$D_{-1}$}}
\put(45.00,21.00){\makebox(0,0)[cc]{$\theta_1$}}
\put(76.00,21.00){\makebox(0,0)[cc]{$\theta_0$}}
\put(105.00,21.00){\makebox(0,0)[cc]{$\theta_{-1}$}}
\put(20.00,30.00){\makebox(0,0)[cc]{0}}
\put(20.00,10.00){\makebox(0,0)[cc]{0}}
\put(100.00,26.00){\vector(0,-1){10.89}}
\put(115.00,35.00){\makebox(0,0)[cc]{$\db$}}
\put(128.00,30.00){\makebox(0,0)[cc]{$\cdots$}}
\put(128.00,10.00){\makebox(0,0)[cc]{$\cdots$}}
\put(30.00,30.00){\vector(-1,0){8.00}}
\put(64.00,30.00){\vector(-1,0){15.93}}
\put(94.00,30.00){\vector(-1,0){17.93}}
\put(124.00,30.00){\vector(-1,0){17.93}}
\put(26.00,10.00){\vector(-1,0){4.00}}
\put(86.00,10.00){\vector(-1,0){4.00}}
\put(57.00,10.00){\vector(-1,0){4.1}}
\put(124.00,10.00){\vector(-1,0){9}}
\end{picture}
\]
is commutative. Similarly as above, a derivation
$\chi \in J_1(n)$ is given by a sequence
$\chi_i: \cb i\to \bod^{n-1,1}(V,\cb{i-1})$, $i\leq 1$, of
$(V,\MU)$-linear maps. The condition $\nabla(\chi)=\theta$ then means
that $\theta_i = D_{-1}\circ \chi_i + \chi_{i+1}\circ \db$,
especially, $\theta_1=D_{-1}\circ \chi_1$. This last equation
implies, since $\chi_1=0$ by~(iii) of Definition~\ref{0}, that
$\nabla(\chi)=\theta$ forces $\theta_1(1_V)=0$. On the other hand, if
$\theta_1(1_V)=0$ then Lemma~\ref{1} gives a $\chi \in J_1(n)$ with
$\theta = \nabla(\chi)$ and we conclude that the contribution of
$\theta_{\bv}$ to $H_0(J_\b(n),\nabla)$ is parametrized by
$\theta_{\bv}(1_V)\in \bod^{n-1,1}(V,\cb 1)_I$ which is the first
summand of~(\ref{2}).
\qed

Let us recall that a (right) {\em differential graded (dg) comp
algebra\/}
(or a {\em nonunital operad\/} in the terminology of~\cite{zebrulka})
is a
bigraded differential space $L =(L_\b(\b),\nabla)$, $L_\b(\b)=
\bigoplus_{k\geq 0,n\geq 2} L_k(n)$,
$\nabla(L_k(n))\subset L_{k-1}(n)$,
together with a system of bilinear operations
\[
\circ_i : L_p(a) \otimes L_q(b) \to L_{p+q}(a+b-1)
\]
given for any $1\leq i \leq b$ such that, for $\phi \in L_p(a)$,
$\psi \in L_q(b)$ and $\nu \in L_r(c)$,
\begin{equation}
\label{4}
\phi\circ_i(\psi\circ_j \nu) =
\left\{
\begin{array}{ll}
(-1)^{p\cdot q}\cdot\psi\circ_{j+a-1}(\phi\circ_i \nu),&
\mbox{ for }1\leq i\leq j-1,
\\
(\phi \circ_{i-j+1} \psi)\circ_j \nu, &
\mbox{ for }j\leq i \leq b+j-1, \mbox{ and}
\\
(-1)^{p\cdot q}\cdot \psi\circ_j(\phi\circ_{i-b+1} \nu),&
\mbox{ for }i\geq j+b.
\end{array}
\right.
\end{equation}
We suppose, moreover, that for any $\phi \in L_p(a)$,
and $\psi \in L_q(b)$,
$1\leq i \leq b$,
\[
\nabla(\phi\circ_t \psi) = \nabla(\phi)\circ_t \psi +
(-1)^p\cdot \phi\circ_t \nabla(\psi).
\]

Any dg comp algebra determines a
nonsymmetric (unital) operad in the monoidal category of differential
graded spaces (see~\cite{may:1972} for the terminology).
To be more precise, let $L = (L_\b(\b),\circ_i,\nabla)$
be a dg comp
algebra as above and let us define the bigraded vector space $
\L_\b(\b) = \bigoplus_{k\geq 0, n\geq 1}$ by $\L_\b(n):=
L_\b(n)$ for $n\geq 2$ and
$\L_\b(1)=\L_0(1):= \mbox{Span}(1_\L)$, where $1_\L$ is a degree zero
generator. Let us extend the definition of structure maps $\circ_i$
to $\L$ by putting $f\circ_1 1_\L:= f$ and $1_\L\circ _i g:= g$, for
$f\in \L_\b(m)$, $g\in \L_\b(n)$ and $1\leq i\leq n$. Let us extend
the differential $\nabla$ by $\nabla(1_\L):=0$. In~\cite{zebrulka} we
proved the following proposition.

\begin{proposition}
The composition maps
$\gamma :\L_\b(a)\ot \L_\b(n_1)\ot \cdots
\ot \L_\b(n_a)\to \L_\b (n_1+\cdots+ n_a)$ given by
\[
\gamma(\phi;\nu_1,\ldots,\nu_a):= \nu_1 \circ_1(\nu_2
\circ_2(\cdots \circ_{a-1}(\nu_a \circ_a \nu)))
\]
for $\phi \in \L_\b(a)$ and $\nu_i \in \L_\b(n_i)$, $1\leq i\leq a$,
define on $\L_\b(\b)$ a structure of a nonsymmetric differential
graded
operad in the monoidal category of differential graded vector spaces.
\end{proposition}

Let $N$ be a (graded) $(V,\MU)$-bimodule, let $d: N\to N$ be a
$(V,\MU)$-linear differential and define
$X_k(n):=\{\theta \in \Der^k_V(\bod(N));\
\theta(N)\subset \bod^n(N)\}$. For $\omega \in X_\b(m)$,
$\theta \in X_\b(n)$ and $1\leq i\leq n$ let $\omega_N
:= \omega|_N : N\to \bod^m(N)$ and $\theta_N := \theta|_N
: N\to \bod^n(N)$ be the restrictions. Let then $\omega \circ_i
\theta \in X_\b(m+n-1)$ be a derivation defined by $(\omega \circ_i
\theta)|_N := (\id^{\odot(i-1)}\odot \omega_N
\odot \id^{\odot (n-i)})\circ \theta_N$. Let us extend the
differential
$d$ to a derivation $D$ of $\bod(N)$ and define $\nabla(\theta):=
[D,\theta]$.

\begin{lemma}
The object $X_\b(\b)=(X_\b(\b), \circ_i,\nabla)$ constructed above
is a differential graded comp algebra.
\end{lemma}

\pf
Let $\phi\in X_p(a)$, $\psi \in X_q(b)$ and $\nu \in X_r(c)$. The
composition $\phi\circ_i(\psi\circ_j\nu)$ is, by definition, given by
its restriction $[\phi\circ_i(\psi\circ_j\nu)]_N$ to $N$ as
\[
[\phi\circ_i(\psi\circ_j\nu)]_N =
\hh{i-1}{\phi}{b+c-i-1}\circ (\psi\circ_j \nu)_N,
\]
with $(\psi\circ_j \nu)_N = \hh{j-1}{\psi}{c-j}\circ \nu_N$.
This implies that
\[
[\phi\circ_i(\psi\circ_j\nu)]_N= \hh{i-1}{\phi}{b+c-i-1}\circ
\hh{j-1}{\psi}{c-j}\circ \nu_N.
\]
For $i\leq j-1$ we have (taking into the account that
$\mbox{\rm Im}(\phi_N)\subset \bod^a(N)$ and
$\mbox{\rm Im}(\psi_N)\subset \bod^b(N)$)
\begin{eqnarray*}
\lefteqn{
\hh{i-1}{\phi}{b+c-i-1}\circ \hh{j-1}{\psi}{c-j}=\hskip.5cm}
\\
&&\hskip.5cm=(-1)^{pq}\cdot \hh{j+a-2}{\psi}{c-j}\circ
\hh{i-1}{\phi}{c-i},
\end{eqnarray*}
which means that $[\phi\circ_i(\psi\circ_j\nu)]_N=
(-1)^{pq}\cdot [\psi\circ_{j+a-1}(\phi\circ_i \nu)]_N$.
This is the axiom~(\ref{4}) for $i\leq j-1$.

Similarly, for $j\leq i\leq b+j-1$ we have
\begin{eqnarray*}
\lefteqn{
\hh{i-1}{\phi}{b+c-i-1}\circ \hh{j-1}{\psi}{c-j}=\hskip.5cm}
\\
&&\hskip.5cm=\id^{\odot (j-1)}\odot
[\hh {i-j}{\phi}{b-i+j-1}\circ \psi_N]\odot
\id^{\odot(c-j)}
\end{eqnarray*}
which means that $[\phi\circ_i(\psi\circ_j\nu)]_N=
[(\phi\circ_{i-j+1}\psi)\circ_j \nu]_N$.
This is the axiom~(\ref{4}) for $j\leq i\leq b+j-1$.
The discussion of the remaining case $i\geq j+b$ is similar.
\qed

Let us consider the special case of the construction above with $N :=
\susp V \oplus \susp \bv$ and $d := 0 \oplus \susp\db\desusp$.

\begin{lemma}
The bigraded subspace $J_\b(\b)$ of $X_\b(\b)$ introduced in
Definition~\ref{0} is closed
under the operations $\circ_i$ and the differential $\nabla$.
\end{lemma}

The proof of the lemma is a straightforward verification.
The lemma says that the dg comp algebra structure
on $X_\b(\b)$ restricts to a dg comp algebra structure
$(J_\b(\b),\circ_i,\nabla)$ on $J_\b(\b)$.

\section{More about $\der *$}

Let $J_\b(\b)= (J_\b(\b),\circ_i,\nabla)$ and $D_{-1}\in \der 1$
be as in the previous section.

\begin{definition}
\label{8}
An infinitesimal deformation of $D_{-1}$ is an element
$D_0\in J_0(2)$ such that $\nabla(D_0) =0$.
An integration of an infinitesimal deformation $D_0$ is a
sequence $\tilde D = \{D_i\in J_i(i+2);\ i\geq 1\}$ such that $D:=
D_{-1}+ D_0 + D_1+ \cdots$ satisfies $[D,D]=0$.
\end{definition}

Let $K_n$ be, for $n\geq 2$, the Stasheff
associahedron~\cite{stasheff:TAMS63}.
It is an $(n-2)$-dimensional cellular complex whose $i$-dimensional
cells are indexed by the set $\Br in$ of all (meaningful) insertions
of $(n-i-2)$ pairs of brackets between $n$ symbols, with suitably
defined incidence maps. There is, for any $a,b \geq 2$, $0\leq i \leq
a-2$, $0\leq j \leq b-2$ and $1\leq t \leq b$, a map
\[
(-,-)_t : \Br ia \times \Br jb \to \Br{i+j}{a+b-1},\ u\times v \mapsto
(u,v)_t,
\]
where $(u,v)_t$ is given by the insertion of $(u)$ at the $t$-th
place in $v$.
This map defines, for $a,b \geq 2$ and $1\leq t\leq b$, the inclusions
$\iota_t :K_a\times K_b \hookrightarrow \partial K_{a+b-1}$.
It is well-known that the sequence $\{K_n\}_{n\geq 1}$ form a
topological operad, see again~\cite{stasheff:TAMS63}.

Let $\cc in$ denote the set of $i$-dimensional oriented
cellular chains with coefficients in $\boldk$ and let $d_C
:\cc in\to \cc {i-1}n$ be the cellular differential.
For $\bolds\in \cc pa$ and $\boldt \in \cc qb$, $p,q \geq 0$,
$a,b\geq 2$ and $1\leq i\leq b$,
let $\bolds \times\boldt \in C_{p+q}(K_a\times K_b)$
denote the cellular cross product and put
\[
\bolds\circ_i\boldt:=
({\iota_i})_*
(\bolds\times\boldt)\in \cc {p+q}{a+b-1}.
\]

\begin{proposition}
\label{5}
The cellular chain complex $(\cc \b\b,d_C)$ together with
operations $\circ_i$ introduced above forms a
differential graded comp algebra.
\end{proposition}

The simplicial version of this proposition was proved
in~\cite{markl-shnider:preprint},
the proof of the cellular version is similar. The dg comp algebra
structure of Proposition~\ref{5} reflects the topological operad
structure of $\{K_n\}_{n\geq 2}$ mentioned above.

Let $c_n$ be, for $n\geq 0$, the unique top dimensional cell of
$K_{n+2}$, i.e.~the unique element of $\Br{n+2}n$ corresponding to
the insertion of no pairs of brackets between $(n+2)$ symbols. Let us
define $\e_n \in \cc n{n+2}$ as $\e_n:= 1\cdot c_n$.

The following proposition was proved in~\cite{zebrulka}, see also the
comments below.

\begin{proposition}
\label{6}
The graded comp algebra $\cc {\b}{\b}= (\cc\b\b,\circ_i)$ is a free
graded
comp algebra on the set $\{\e_0,\e_1,\ldots\}$.
\end{proposition}

Let us recall that the freeness in the proposition above means that
for any graded comp algebra $L_\b(\b)= (L_\b(\b),\circ_i)$ and
for any sequence $\alpha_n\in L_n(n+2)$, $n\geq 0$, there exists a
unique graded comp algebra map $f : \cc {\b}{\b}\to L_\b(\b)$ such
that $f(\e_n)= \alpha_n$, $i\geq 0$.

The proof of Proposition~\ref{6}
is based on the following observation. There is a description of the
free graded comp algebra (= free nonsymmetric
nonunital operad) on a given set in terms of oriented planar trees.
The free comp algebra ${\cal F}(\e_0,\e_1,\ldots)$ on
the set $\{\e_0,\e_1,\ldots\}$ has
${\cal F}(\e_0,\e_1,\ldots)(n) =$ the vector space spanned by oriented
connected planar trees with $n$ input edges. Each such a tree $T$ then
determines an element of $\Br ni$ where $i =$ the number of
vertices of $T$, i.e.~a cell of $\cc{n-i-2}n$. This correspondence
defines a map ${\cal F}(\e_0,\e_1,\ldots)(n) \to \cc *n$ which
induces the requisite isomorphism of operads.

Let $L= (L_\b(\b),\circ_i,\nabla)$ be a dg comp
algebra. Let us define, for $\phi \in L_p(a)$ and $\psi \in L_q(b)$,
\[
\phi\diamond \psi :=
\sum_{1\leq i\leq b}(-1)^{(a+1)(i+q+1)}\cdot \phi\circ_i \psi
\mbox{\hskip2mm and \hskip2mm}
[\phi,\psi]:= \phi\diamond \psi -
(-1)^{(a+p+1)(b+q+1?)} \cdot\psi\diamond \phi.
\]
In~\cite{markl-shnider:preprint} we proved the following proposition.

\begin{proposition}
The operation $[-,-]$ introduced above endows $L^\b:=
\bigoplus_{a-p-1=\b}L_p(a)$ with a structure of a differential graded
(dg) Lie algebra, $L = (L^\b,[-,-],\nabla)$.
\end{proposition}

The construction above thus defines a functor from the category of
dg comp algebras to the category of dg Lie algebras.
Let us observe that for the comp algebras $X_\b(\b)$
and $J_\b(\b)$ this structure coincides with the Lie algebra
structure induced by the graded commutator of derivations.
We can also easily prove that the elements $\{\e_n\}_{n\geq 0}$
satisfy, for each $m\geq 0$,
\begin{equation}
\label{7}
\mbox{$
d_C(\e_m)+\frac12 \sum_{i+j=n-1}[\e_i,\e_j]=0
$}
\end{equation}
in the dg Lie algebra $CC(K)^*=(CC(K)^\b,[-,-],d_C)$.

\begin{proposition}
There exists an one-to-one correspondence between integrations of an
infinitesimal deformation $D_0$ in the sense of
Definition~\ref{8} and dg comp algebra homomorphisms
$m: \cc \b\b \to J_\b(\b)$ with $m(\e_0)=D_0$.
\end{proposition}

\pf
Let us suppose we have a map $m: \cc \b\b
\to J_\b(\b)$ with $m(\e_0)=D_0$
and define, for $n\geq 1$, $D_n:= m(\e_n)$.
We must prove that the derivation $D:= D_{-1}+ D_0 + D_1 +\cdots$
satisfies $[D,D]=0$. This condition means that
\begin{equation}
\label{9}
\mbox{$
\nabla(D_m)+\frac12 \sum_{i+j=n-1}[D_i,D_j]=0
$}
\end{equation}
for any $m\geq 0$, which is exactly what we get
applying on~(\ref{7}) the Lie algebra homomorphism $m$.

On the other hand, suppose we have an integration
$\{D_n\}_{n\geq 1}$. The freeness of the graded comp algebra $\cc
\b\b$
(Proposition~\ref{6}) ensures the existence of a graded comp algebra
map $m: \cc \b\b \to J_\b(\b)$ with $m(\e_n)= D_n$ for $n\geq 0$. We
must verify that this unique map commutes with the differentials,
i.e.~that $m(d_C(\bolds))= \nabla(m(\bolds))$ for any $\bolds\in
\cc**$.
Because of the freeness, it is enough to verify the last condition for
$\bolds\in \{\e_n\}_{n\geq 0}$, i.e.~to verify that $m(d_C(\e_n))=
\nabla(m(\e_n))$ for $n\geq 0$. Expanding $d_C(\e_n)$
using~(\ref{7}) we see that this follows from~(\ref{9}) and from the
fact that the map $m$ is a homomorphism of graded Lie algebras.
\qed

\begin{proposition}
\label{13}
An infinitesimal deformation $D_0\in J_0(2)$ can be
integrated if and only if
$[D_0,D_0]|_{\susp V} =[D_0,D_0]|_{\susp \calb_1(V)}=0$.
\end{proposition}

\pf
Standard obstruction theory.
Let us suppose that we have an integration $\tilde D =
\{D_i\}_{i\geq 1}$. Condition~(\ref{9}) with $m=1$ means that
\begin{equation}
\label{11}
\mbox{$\nabla(D_1)+\frac12[D_0,D_0]=0.$}
\end{equation}
Because $D_1\in J_1(3)$,
$D_1|_{\susp V} = D_1|_{\susp \calb_1(V)}=0$ by (iii) of
Definition~\ref{0},
therefore $\nabla(D_1)|_{\susp V} =\nabla(D_1)|_{\susp \calb_1(V)}=0$
and~(\ref{11}) implies that $[D_0,D_0]|_{\susp V} =
[D_0,D_0]|_{\susp \calb_1(V)}=0$.

On the other hand, let us suppose that $[D_0,D_0]|_{\susp V}
=[D_0,D_0]|_{\susp \calb_1(V)}=0$. By the description of
$H_0(J_*(3),\nabla)$ as it is given in Proposition~\ref{12} we see
that the homology class of $[D_0,D_0] \in J_0(3)$ is zero, therefore
there exists some $D_1\in J_1(3)$ with
$\nabla(D_1)+\frac12[D_0,D_0]=0$.

Let us suppose that we have already constructed a sequence
$D_i\in J_i(i+2)$, $1\leq i\leq N$, such that basic equation~(\ref{9})
holds for any $m\leq N$. The element $\frac12
\sum_{i+j=N}[D_i,D_j]\in J_N(N+3)$ is a $\nabla$-cycle (this follows
from the definition of $\nabla(-)$ as $[D_{-1},-]$ and the Jacobi
identity) and the triviality of $H_N(J_*(N+3))$ (again
Proposition~\ref{12}) gives some $D_{N+1}\in J_{N+1}(N+3)$ which
satisfies~(\ref{9}) for $m=N+1$. The induction may go on.
\qed

Let $G$ be the subgroup of $\mbox{\rm Aut}(\bod(V,\bv)^*)$ consisting
of automorphisms of the form $g= \id+\phi_{\geq 2}$, where
$\phi_{\geq 2}$ is a $(V,\MU)$-linear map satisfying
$\phi_{\geq 2}(\bod^{i,j}(V,\bv))\subset\bod^{\geq i+2,j}(V,\bv)$ for
any $i,j$.

Let us observe that $G$ naturally acts on the set of integrations of
a fixed infinitesimal deformation $D_0$. To see this, let $\tilde D =
\{D_i\}_{i\geq 1}$ be such an integration and let us denote, as
usual, $D := D_{-1}+D_0+D_1 +\cdots$. Then $g^{-1}Dg$ is, for
$g\in G$, clearly also a $(V,\MU)$-linear derivation from $\der 1$
and $(g^{-1}Dg)(\bod^{i,j}(V,\bv))\subset (\bod^{\geq i,j}(V,\bv))$.
We may thus decompose $g^{-1}Dg$ as $g^{-1}Dg = \sum_{k\geq -1}D'_k$
with $D'_k(\bod^{i,j}(V,\bv))\subset \bod^{i+k+1,j}(V,\bv)$. We
observe that $D'_{-1}=D_{-1}$, $D'_0=D_0$ and that, from degree
reasons, $D'_k(\bod(V,\bv)^{i,j})\subset \bod(V,\bv)^{i+k+1,j-k}$.
This means that $D'_k \in J_k(k+2)$ for $k\geq 1$. The equation
$[g^{-1}Dg,g^{-1}Dg]=0$ is immediate, therefore the correspondence
$(g,\{D_i\}_{i\geq 1}) \mapsto \{D'_i\}_{i\geq 1}$ defines the
requisite action. The following proposition shows that this action is
transitive.

\begin{proposition}
\label{16}
Let $\tilde D' = \{D'_i\}_{i\geq 1}$ and $\tilde D'' =
\{D''_i\}_{i\geq 1}$ be two integrations of an infinitesimal
deformation $D_0$. If we denote $D' := D_{-1}+D_0+\sum_{i\geq 1}D'_i$
and $D'' := D_{-1}+D_0+\sum_{i\geq 1}D''_i$, then $D' = g^{-1}D''g$
for some $g\in G$.
\end{proposition}

\pf
Again standard obstruction theory.
As we already observed, for any $g\in G$, $g^{-1}D''g$ decomposes
as $g^{-1}D''g = D_{-1}+D_0+\sum_{i\geq 1}\{g^{-1}D''g\}_i$ with some
$\{g^{-1}D''g\}_i \in J_i(i+2)$. Let us suppose that we have already
constructed some $g_N \in G$, $N\geq 1$, such that
$\{g_N^{-1}D''g_N\}_i= D'_i$ for $1\leq i \leq N$. We have
\[
\mbox{$\nabla(D'_{N+1})-\frac12 \sum_{i+j=N}[D'_i,D'_j]=0$}
\]
and, similarly,
\[
\mbox{$\nabla(\{g_N^{-1}D''g_N\}_{N+1})-\frac12
\sum_{i+j=N}[\{g_N^{-1}D''g_N\}_i,\{g_N^{-1}D''g_N\}_j]=0$.}
\]
By the induction, the second terms of the above equations are the
same, therefore $\nabla(D'_{N+1}) = \nabla(\{g_N^{-1}D''g_N\}_{N+1})$
which means that $D'_{N+1} -
\{g_N^{-1}D''g_N\}_{N+1}\in J_{N+1}(N+3)$ is a cycle. The triviality
of $H_{N+1}(J_*(N+3))$ (Proposition~\ref{12}) gives some
$\phi \in J_{N+2}(N+3)$ such that $D'_{N+1} -
\{g_N^{-1}D''g_N\}_{N+1}= \nabla(\phi)$.
The element $\mbox{\rm exp}(\phi)\in G$ is of the form $\id +
\phi +\phi_{\geq N+3}$ with
\[
\phi_{\geq N+3}(\bod^{i,j}(V,\bv))\subset
\bod^{i+N+3,j}(V,\bv),
\]
therefore $g_{N+1}:= \mbox{\rm exp}(\phi)g_N$ satisfies
$\{g_{N+1}^{-1}D''g_{N+1}\}_i= \{g_{N}^{-1}D''g_{N}\}_i
=D'_i$ for $1\leq i \leq N$ and $\{g_{N+1}^{-1}D''g_{N+1}\}_{N+1}=
\{g_{N}^{-1}D''g_{N}\}_{N+1}+\nabla(\phi)= D'_{N+1}$ and the
induction goes on. The prounipotency of the group $G$ assures that
the sequence $\{g_N\}_{N\geq 1}$ converges to some $g\in G$ as
required.
\qed

\section{Applications to Drinfel'd algebras}

In~\cite{markl-shnider:preprint} we
introduced two $(V,\MU)$-linear `coactions' $\lambda
:\bv \to V\odot \bv$ and $\rho : \bv \to \bv\odot V$ as
\begin{eqnarray*}
\lambda(a_0|\cdots|a_{-n+1})&:=&
\sum\Delta'(a_0)\cdots\Delta'(a_{-n+1})\od
(\Delta''(a_0)|\cdots|\Delta''(a_{-n+1})),\mbox{ and}
\\
\rho(a_0|\cdots|a_{-n+1})&:=&\sum
(\Delta'(a_0)|\cdots|\Delta'(a_{-n+1}))\od\Delta''(a_0)\cdots
\Delta''(a_{-n+1}).
\end{eqnarray*}
Let us define a derivation $D_0\in J_0(2)$ by
\[
D_0|_{\susp \bv}:= (\susp\bod\susp)(\lambda+\rho)(\desusp)
\mbox{ and }
D_0|_{\susp V}:= (\susp\bod\susp)(\Delta)(\desusp).
\]

\begin{proposition}
The derivation $D_0$ defined above is an integrable infinitesimal
deformation of $D_{-1}$.
\end{proposition}

\pf
To prove that $D_0$ is an infinitesimal deformation of $D_{-1}$ means
to show that $\nabla(D_0)= [D_{-1},D_0]= 0$.
This was done in~\cite{markl-shnider:preprint}.

The integrability of $D_0$ means, by Proposition~\ref{13}, that
$[D_0,D_0]|_{\susp V}= [D_0,D_0]|_{\susp \calb_1(V)}=0$. For $v\in V$
we have
\begin{eqnarray*}
[D_0,D_0](\susp v)&=& (D_0 \circ D_0)(\susp v)= D_0 (\susp
\odot \susp)(\Delta)(\susp v)=
\\
&=&[((\susp \odot \susp)(\Delta)\odot \susp)(\Delta)-
(\susp \odot(\susp \odot \susp)(\Delta))(\Delta)](\susp v) =
\\
&=&(\susp\odot\susp\odot\susp)\circ [(\Delta\odot \id)(\Delta)-
(\id \odot \Delta)(\Delta)(\susp v)]
\end{eqnarray*}
which is zero by the quasi-coassociativity~(\ref{14}) and
by~(\ref{3}).

Similarly, for $(v)\in \calb_1(V)$ we have
\[
D_0^2 = (\susp \odot \susp \odot \susp)[ (\Delta\odot \id)\lambda -
(\id \odot \lambda)\lambda -(\id \odot\rho)\lambda
+ (\lambda \odot \id)\rho + (\rho \odot \id)\rho
-(\id \odot \Delta)\rho ](\desusp (v))
\]
and~(\ref{14}),~(\ref{3}) again imply that this is zero.
\qed

\begin{definition}
\label{15}
An integration of the infinitesimal deformation $D_0$ above is called
a homotopy comodule structure.
\end{definition}

Let $\bod'(V,\bv)= \bod^{*,1}(V,\bv)$ denote the submodule of
$\bod^{*,*}(V,\bv)$ with precisely one factor of $\susp \bv$.
Let $C^n(A)$ be the set of all degree $n$ homogeneous maps $f:
\bod'(V,\bv) \to \bod^*(\susp V)$ which are both $\bod^*(\susp V)$
and $(V,\MU)$-linear. Let us define also a degree one derivation
$d_C$ on $\bod^*(\susp V)$ by $d_C|_{\susp V}:=
(\susp \odot \susp)(\Delta)(\susp)$.

Let $\{D_i\}_{i\geq 1}$ be a homotopy comodule structure in the sense
of Definition~\ref{15} and let $D:= D_{-1}+D_0+D_1+\cdots$. Define a
degree one endomorphism $d$ of $C^*(A)$ by $d(f):= f\circ D+ (-1)^n
d_C\circ f$. It is easy to show that $d$ is a differential and,
following~\cite{markl-shnider:preprint}, we define the
{\em cohomology of our Drinfel'd algebra
$A$\/} by $H^*(A):= H^*(C^*(A),d)$.

\begin{proposition}
The definition of the cohomology of a Drinfel'd algebra does not
depend on the particular choice of a homotopy comodule structure.
\end{proposition}

\pf
Let $\{D'_i\}_{i\geq 1}$ and $\{D''_i\}_{i\geq 1}$ be two homotopy
comodule structures, $D':= D_{-1}+D_0+D'_1+\cdots$ and
$D'':= D_{-1}+D_0+D''_1+\cdots$. Let $d'(f):= f\circ D'+ (-1)^n
d_C\circ f$ and $d'(f):= f\circ D''+ (-1)^n
d_C\circ f$. Proposition~\ref{16} then gives
some $g\in G$ such that $D' \circ g = g \circ D''$. We see
immediately that the map $\Psi :(C^*(A),d') \to (C^*(A),d'')$ defined
by $\Psi(f):= f\circ g$ is an isomorphism of complexes.
\qed

\vskip3mm
\catcode`\@=11
\noindent
M.~M.: Mathematical Institute of the Academy, \v Zitn\'a 25, 115 67
Praha 1, Czech Republic,\hfill\break\noindent
\hphantom{M.~M.:}\hskip1mm
email: {\bf markl@earn.cvut.cz}\hfill\break\noindent
\noindent

\noindent
S.~S.: Department of Mathematics, University of Bar Ilan, Israel,
\hfill\break\noindent
\hphantom{S.~S.:}\hskip1mm email: {\bf shnider@bimacs.cs.biu.ac.il}

\end{document}